\documentclass[12pt]{article}

\usepackage{psfig}

\usepackage{amssymb}
\def\suffix{ps}
\newcount\system
\global\system=3   % for unix(dvips)
\def\ifundefined#1{\expandafter\ifx\csname#1\endcsname\relax}
\ifundefined{figdir}\def\figdir{}\fi
\newcount\firstline
\newdimen\pswidth  \newdimen\xleft
\newdimen\psheight \newdimen\ytop \newdimen\ybot
\newcount\justx \newcount\justy
\global\justx=0 \global\justy=0
\newdimen\vpos \newtoks\labeL 
\newread\labeLfile \newdimen\xcoord \newdimen\ycoord
\newif\ifdoit 
\newbox\labox
\newdimen\xdvikwid 
\newdimen\xdvikht
\newdimen\pspoints
\newdimen\rwi
\newcount\temp
\def\readdim#1{\global\read\labeLfile to \temp
\global #1=\temp pt}
\def\figcrop#1{\par%  #1=filename
\openin\labeLfile=\figdir#1.lbl                                              
\global\read\labeLfile to\firstline\message{#1}               
\global\read\labeLfile to\temp%read overall dimensions                                     
\readdim{\ybot}
\readdim{\xleft}%               read upper left point
\readdim{\ytop}
\global\read\labeLfile to\justx%ignore
\global\read\labeLfile to\justy%ignore
\global\read\labeLfile to\labeL%ignore
\readdim{\pswidth}%            read lower right point
\global\advance\pswidth by -\xleft
\readdim{\psheight}
\global\advance\ybot by -\psheight
\global\advance\psheight by -\ytop
\global\read\labeLfile to\justx%ignore
\global\read\labeLfile to\justy%ignore
\global\read\labeLfile to\labeL%ignore                                    
\vbox to\psheight{\vfill
\ifnum\system=1
\fi %textures
\ifnum\system=2
\fi %msdos
\ifnum\system=3
                                                 \fi         %%unix:dvips
\ifnum\system=4
\fi         %%unix:dvips,scaled
\ifnum\system=1
\hbox to \pswidth{\kern-\xleft\special{postscriptfile \figdir#1.\suffix }\hfil}\fi
\ifnum\system=2
\hbox to \pswidth{\kern-\xleft\special{ps: plotfile \figdir#1.\suffix }\hfil}\fi
\ifnum\system=3
\hbox to \pswidth{\kern-\xleft\includegraphics{\figdir#1.\suffix}\hfil}\fi
\ifnum\system=4
\hbox to \pswidth{\kern-\xleft\includegraphics{\figdir#1.\suffix}\hfil}\fi
\ifnum\system=5
\hbox to \pswidth{\kern-\xleft\includegraphics{\figdir#1.\suffix}\hfil}\fi %orphee
\ifnum\system=6
   \xdvikwid=\pswidth
   \xdvikht=\psheight
   {\global\divide\xdvikwid by \pspoints}
   {\global\divide\xdvikht by \pspoints}
   \rwi=\xdvikwid
    {\global\multiply\rwi by 10}
\hbox to \pswidth{\kern-\xleft\includegraphics{\figdir#1.\suffix\space}\hfil}\fi                   %xdvik
\vskip -\baselineskip
\vskip -\ybot 
\vskip-\psheight %                                     
\hbox to\pswidth  {\hss}%                                            
\parindent=0pt\offinterlineskip                                       
\vpos=0 pt%                                                              
\loop\readdim{\xcoord}                                 
\ifdim \xcoord < -999pt \doitfalse\else\doittrue\fi                        
\ifdoit \advance \xcoord by -\xleft
\readdim{\ycoord}
\advance \ycoord by -\ytop                              
\global\read\labeLfile to\justx                                       
\global\read\labeLfile to\justy                                       
\global\read\labeLfile to\labeL
\global\setbox\labox=\hbox{\labeL\hskip-0.3em}%    
\advance\vpos by-\ycoord                                              
\vskip-\vpos \vpos=\ycoord                                         
\hbox to\pswidth{\hskip\xcoord %                                 
\hbox to 0pt{\ifnum\justx>0\hss\fi%                                   
\vbox to0pt{%                                                         
\ifnum\justy<2\vss\fi%                                                
\copy\labox\kern0pt%  
\ifnum\justy>0\vss\fi}%                                               
\ifnum\justx<2\hss\fi}%                                               
\hss}%                                                                
\repeat%                                                              
\advance\vpos by-\psheight%                                           
\vskip-\vpos %                                                     
}\closein\labeLfile}
\def\figplace#1#2#3{
\openin\labeLfile=\figdir#1.lbl
\ifeof \labeLfile
       \immediate\write16{***Can't find \figdir#1.lbl; Skipping it.***}
\else  \closein\labeLfile
       \null\hskip#2\raise #3 \hbox{\figcrop{#1}}
\fi
}
\def\figput#1{
\openin\labeLfile=\figdir#1.lbl
\ifeof \labeLfile
       \immediate\write16{***Can't find \figdir#1.lbl; Skipping it.***}
\else  \closein\labeLfile
       \hbox{\figcrop{#1}}
\fi
}

\def\CC{\mathbb{C}}

\newfam\gotfam
\font\twlgot=eufm10 at 12pt \font\tengot=eufm10
 \font\sevengot=eufm7
\textfont\gotfam=\twlgot \scriptfont\gotfam=\tengot
\scriptscriptfont\gotfam=\sevengot
\def\got{\fam\gotfam\twlgot}

\newtheorem{theor}   {Theorem}

\newcommand{\be}  {\begin{equation}}
\newcommand{\ee}  {\end{equation}}
\newcommand{\bea} {\begin{eqnarray}}
\newcommand{\eea} {\end{eqnarray}}
\newcommand{\lp}  {\left(}
\newcommand{\rp}  {\right)}

\newcommand{\Br}  {\overline}

\newcommand{\cT}  {{\cal T}}

\newcommand{\cA}  {{\cal A}}

\newcommand{\cZ}  {{\cal Z}}

\newcommand{\Om}  {\Omega}

\newcommand{\Ga}  {\Gamma}

\newcommand{\al}  {\alpha}

\newcommand{\de}  {\delta}

\newcommand{\ph}  {\phi}

\newcommand{\Gc}  {{\got c}}

\newcommand{\Gl}  {{\got l}}

\def\Br{\overline}

\newcommand{\eqdef} {\stackrel{\rm def}{=}}

\def\endproof{\hfill\vrule height .6em width .6em depth
  0pt\goodbreak\vskip.25in}

\begin{document}

\title{The Jacobian Conjecture as a Problem of Perturbative Quantum
Field Theory}

\author{Abdelmalek Abdesselam \\
\\
{\small D{\'e}partement de Math{\'e}matiques}\\
{\small  Universit{\'e} Paris XIII, Villetaneuse}\\
{\small Avenue J.B. Cl{\'e}ment, F93430 Villetaneuse, France}\\
}

\maketitle

%\hfill\eject

{\abstract{
The Jacobian conjecture is an old unsolved
problem in mathematics, which has been unsuccessfully
attacked from many different angles. We add here another
point of view pertaining to the so called formal inverse approach,
that of perturbative quantum field theory.} }

 
%\hfill\eject

\section{Introduction}

The purpose of this modest note, for which we claim no originality
except that of connecting apparently unrelated fields,
is to draw the attention of theoretical physicists to one of the major
unsolved problems of mathematics~\cite{Smale}, viz. the Jacobian conjecture.
The question is so simple that it was coined in~\cite{Abhyankar1}
a problem in ``high school algebra''. One can formulate it as follows.

Let $F:\CC^n\rightarrow \CC^n$ be a map written in coordinates as
\be
F(x_1,\ldots,x_n)=(F_1(x_1,\ldots,x_n),\ldots,F_n(x_1,\ldots,x_n))
\ .
\ee
One says that $F$ is a {\em polynomial map} if the functions
$F_i:\CC^n\rightarrow \CC$ are polynomial.
Suppose that the Jacobian determinant
\be
JF(x_1,\ldots,x_n)\eqdef
det\lp \frac{\partial F_i}{\partial x_j}(x_1,\ldots,x_n)\rp
\ee
is {\em identically} equal to a nonzero constant. Show then that $F$
is globally invertible (for the composition of maps) and that its inverse
$G\eqdef F^{-1}$ is also a polynomial map.

Since it was first proposed in~\cite{Keller} (for $n=2$ and polynomials with
integral coefficients), this problem has resisted all attempts for a solution.
In fact, this seemingly simple problem is quite an embarrassment.
Indeed, some faulty proofs have even been published (see the
indispensable~\cite{BassCW} and~\cite{vdEssen} for a review).
We will show here that the Jacobian conjecture
can be formulated in very nice way as a 
question in preturbative quantum field theory (QFT).
We also expect any future progress on this question to be beneficial
not only for mathematics, but also for theoretical physics as it
would enhance our understanding of perturbation theory. 

\bigskip
\noindent{\bf Aknowledgments :}
The author is grateful to V. Rivasseau for early encouragements and
collaboration on this project. Some of the ideas presented here
are due to him.
We thank D. Brydges, C. de Calan and J. Magnen for enlighting discussions.
We also thank J. Feldman for his invitation to the Mathematics Department
of the University of British Columbia where part of this work was done.
The pictures in this article have been drawn using a software
package kindly provided by J. Feldman.

\section{The formal inverse as a one-point correlation function}

The most tempting, yet unfortunately least developed, line
of attack on the Jacobian conjecture is the so called formal inverse
approach.
One tries to solve explicitly for $x=(x_1,\ldots,x_n)$
in the equation $y=F(x)$, one then finds a power series expression for $x$
in terms of $y=(y_1,\ldots,y_n)$. By the uniqueness of the power series
inverse, all one has to do then is to show that it is in fact a polynomial,
that is the terms of high degree in the $y$ variables vanish.
One of the many reasons this approach is in its infancy is that it took
more than two centuries (say from~\cite{Lagrange} to~\cite{BassCW})
to have a workable
formula for the formal inverse in the multivariable case.
Early contributions can be found
in~\cite{Laplace,Jacobi,Darboux,Stieltjes,Poincare,Good}.
An important contribution concerning formal inversion
is due to Gurjar and Abhyankar~\cite{Abhyankar2}. 
Modern litterature on reversion and Lagrange-Good type formulas
is huge and we invite the reader to
consult~\cite{BassCW,Gessel,Henrici,Wright1} 
for more complete references.
The first formula for the coefficients of the formal inverse power series $G$
in terms of those of $F$ is due to J. Towber and was first published
in~\cite{Wright2}. In physicists' terms ours is the following.

\medskip
\noindent{\bf Claim : (A. A., V. Rivasseau)}
The formal solution of $y=F(x)$, without any assumption on $F$
except that its linear part is invertible,
is the pertubation expansion of the normalized one-point correlation
function
\be
x_i=\frac{1}{Z}\int_{\CC^n}d{\Br \ph}d\ph\ 
\ph_i e^{-{\Br \ph}F(\ph)+{\Br\ph}y}
\label{F1}
\ee
where ${\Br\ph}_1,\ldots{\Br\ph}_n$,$\ph_1,\ldots,\ph_n$
are the components of a complex Bosonic field.
The integration is over $\CC^n$ with the measure
\be
d{\Br \ph}d\ph
\eqdef
\prod_{i=1}^n
\lp
\frac{d(Re\ \ph_i)d(Im\ \ph_i)}{\pi}
\rp
\ ,
\ee
we used the notation ${\Br \ph}F(\ph)\eqdef\sum_{i=1}^n
{\Br\ph}_i F_i(\ph_1,\ldots,\ph_n)$, ${\Br\ph}y\eqdef
\sum_{i=1}^n {\Br\ph}_i y_i$, and
\be
Z\eqdef
\int_{\CC^n}d{\Br \ph}d\ph\ 
e^{-{\Br \ph}F(\ph)+{\Br\ph}y}
\ee

\medskip
We obtained this expression by solving iteratively the equation
$y=F(x)$ thereby generating a tree expansion in the same way
one expresses the effective action $\Ga(\ph)$ in terms of the logarithm
$W(J)$ of the partition function in QFT (see~\cite{Zinn} for instance).
We then determined the Feynman rules of this tree expansion and finally the
``path integral'' formulation (\ref{F1}), only to realize that in fact
our formula is closely related to the one introduced
by G. Gallavotti, following
a suggestion of G. Parisi, to express the Lindstedt perturbation series in
the context of KAM theory~\cite{Gallavotti}.

A mathematician will undoubtedly shriek at the sight of equation
(\ref{F1}). In the following, we will state and prove a precise theorem,
using some analysis, for the case where $F_i(x)=x_i-H_i(x)$, with
the $H_i(x)$, $1\le i\le n$, being homogenous of the same degree $d$.
Indeed, it is enough to treat the cubic case $d=3$, {\em for all dimensions}
$n$, in order to prove the Jacobian conjecture in full
generality~\cite{BassCW}.
However, formula (\ref{F1}) is completely combinatorial in nature and its
proper setting is in the ring of formal power series with variables
corresponding to the coefficients of $F$ together with the $y_i$'s,
over any field of zero characteristic.
One simply has to define {\em formal Gaussian integration}, somewhat
in the spirit of~\cite{BarNatan1,BarNatan2}.
We refer to the expository article~\cite{Abdesselam1}
for a formulation and proof of our claim as a decent mathematical theorem.
The latter article will also provide more details 
on how Feynman diagrams can be useful
in algebraic combinatorics and how well they fit in the Joyal theory
of combinatorial species~\cite{Joyal}. We also refer
to~\cite{Abdesselam2}
for a very simple heuristic proof of the Lagrange-Good multivariable
inversion formula, which becomes a fully rigorous and purely combinatorial
proof when interpreted using the formalism of~\cite{Abdesselam1}.

Now let $F_i(x)=x_i-H_i(x)$ with $H_i(x)$ written in {\em tensorial}
notation as
\be
H_i(x)=\frac{1}{d!}
\sum_{j_1,\ldots,j_d=1}^n
w_{i,j_1\ldots j_d} x_{j_1}\ldots x_{j_d}
\ee
so that the $1$-contravariant and $d$-covariant tensor
$w_{i,j_1\ldots j_d}$ is {\em completely symmetric}
in the $j$ indices.
Let us write
\be
{\Br\ph}w\ph^d\eqdef
\sum_{i=1}^n\sum_{j_1,\ldots,j_d=1}^n
{\Br\ph}_i
w_{i,j_1\ldots j_d} \ph_{j_1}\ldots \ph_{j_d}
\ee
so that (\ref{F1})
becomes
\be
G_i(y)=\frac{\int d{\Br\ph}d\ph\ \ph_i
\exp\lp -{\Br\ph}\ph+\frac{1}{d!}{\Br\ph}w\ph^d
+{\Br\ph}y\rp}
{\int d{\Br\ph}d\ph\ 
\exp\lp -{\Br\ph}\ph+\frac{1}{d!}{\Br\ph}w\ph^d
+{\Br\ph}y\rp}
\label{F2}
\ee
The free propagator is represented as an oriented line
\be
\figplace{propag}{0 in}{-0.06 in}=\de_{ij}
\label{Fig1}
\ee
for the contraction of a pair $\ph_i{\Br\ph}_j$.
There are two types of vertices: the $w$-vertices represented
by
\be
\figplace{wvertex}{0 in}{-0.45 in} d {\rm\ half\ lines}={\Br\ph}w\ph^d
\label{Fig2}
\ee
and the $y$-vertices represented by
\be
\figplace{yvertex}{0 in}{-0.47 in}={\Br\ph}y\ .
\ee
As is well known in QFT, the numerator and denominator of (\ref{F2})
can be calculated by expanding
\[
\exp\lp \frac{1}{d!}{\Br\ph}w\ph^d
+{\Br\ph}y\rp
\]
and integrating term by term with respect to the normalized complex
Gaussian measure $d{\Br\ph}d\ph\ e^{-{\Br\ph}\ph}$.
The result is a sum over all possible Feynman diagrams
that can be built from the vertices of (\ref{Fig1}) and
(\ref{Fig2})
(and the source $\ph_i$ for the numerator) by joining the half-lines
of compatible directions.
A quick look at the vertices shows that the only possible diagrams
are trees connected to the source $\ph_i$, or vacuum graphs made
by an oriented loop of say $k\ge 1$ $w$-vertices to which
$k(d-1)$ trees, whose leaves are $y$-vertices, are attached.
When one factors out the denominator, the only diagrams that remain are
made of a single tree with the source $\ph_i$ as its root.
Therefore, at least formally, we have
\be
G_i(y)=\sum_{V\ge 0}\sum_{N\ge 0}
\frac{1}{V!N!}\sum_{\cT}
\cA_i(\cT)
\label{F5}
\ee
where $\cT$ is a Cayley tree (viewed as a set of unordered pairs)
on a finite set $E=E(V,N)$. The latter is chosen, non canonically, once
for each pair $(V,N)$, and must be the disjoint union of $E_{\rm root}$
of cardinality $1$, $E_{\rm internal}$ of cardinality $V$ and
$E_{\rm leaf}$ of cardinality $N$.
$\cT$ is constrained by the condition that elements of
$E_{\rm root}\cup E_{\rm leaf}$
have valence $1$ while those of $E_{\rm internal}$ have valence $d+1$.
This automatically enforces the relation $(d-1)V=N-1$ which
can be checked by counting the half-lines.
Even though we write, in the sequel, seemingly independent sums over
$V$ and $N$, the previous relation is allways assumed.

We now define the {\em amplitude} $\cA_i(\cT)$.
One directs the edges of $\cT$ towards the root in $E_{\rm root}$.
For each such edge $l\in\cT$, one introduces an index
$i_l$ in the set $\{1,\ldots,n\}$.
One then considers the expression $A_i(\cT,(i_l)_{l\in\cT})$
which is the product of the following factors.

- For each $a\in E_{\rm leaf}$, if $l(a)$ is the unique line going
from $a$, we take the factor $y_{i_{l(a)}}$.

- For each $b\in E_{\rm internal}$, if $\{l_1(b),\ldots,l_d(b)\}$
is the set of lines coming into $b$ and $l_0(b)$ is the unique line
leaving $b$, we take the factor
$w_{i_{l_0(b)},i_{l_1(b)}\ldots i_{l_d(b)}}$.

The resulting monomial in the $y$'s and $w$'s is $A_i(\cT,(i_l)_{l\in\cT})$
by definition.
Finally $\cA_i(\cT)$
is the sum of $A_i(\cT,(i_l)_{l\in\cT})$ over all the indices
$(i_l)_{l\in\cT}$ except the index of the line arriving at the root which
is fixed at the value $i$, the source index.

For example, with $d=3$, the amplitude of the following tree with $V=4$ and
$N=9$
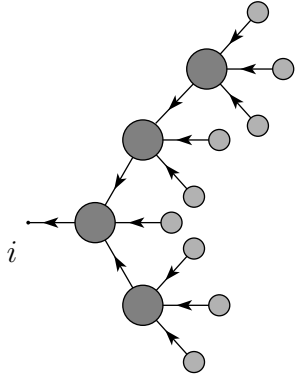
\[
\figput{tree}
\]
is
\bea
\lefteqn{
\cA_i(\cT)=
\sum_{\al_1,\ldots,\al_{12}=1}^n
w_{i,\al_{1}\al_{2}\al_{3}}
w_{\al_{1},\al_{4}\al_{5}\al_{6}}
w_{\al_{3},\al_{7}\al_{8}\al_{9}}
} & & \nonumber\\
 & &
w_{\al_{9},\al_{10}\al_{11}\al_{12}}
y_{\al_2}y_{\al_4}y_{\al_5}
y_{\al_6}y_{\al_7}y_{\al_8}
y_{\al_{10}}y_{\al_{11}}y_{\al_{12}}\ .
\eea

Note that, by the Cayley formula for the number of trees with
preassigned valences, the sum over $\cT$ has
\be
\frac{((1+V+N)-2)!}{(1-1)!((d+1)-1)!^V (1-1)!^N}
=\frac{(V+N-1)!}{d!^V}\ \ {\rm terms}.
\ee
Let us introduce the norms
\be
||w||_{\infty,1}\eqdef
\max_{1\le i\le n}
\sum_{j_1,\ldots,j_d=1}^n
|w_{i,j_1\ldots j_d}|
\ee
and
\be
||y||_{\infty}\eqdef
\max_{1\le i\le n} |y_i|\ .
\ee
We now have
\begin{theor}
The series
\be
G_i(y)=\sum_{V\ge 0}\sum_{N\ge 0}
\frac{1}{V!N!}\sum_{\cT}
\cA_i(\cT)
\ee
is absolutely convergent, provided
\be
||y||_{\infty}< R\eqdef {\lp
\frac{d!}{2^d ||w||_{\infty,1}}
\rp}^{d-1}
\ee
and satisfies, on this domain of convergence,
\be
||G(y)||_{\infty}\le
\frac{||y||_{\infty}}{1-\frac{||y||_{\infty}}{R}}
\label{F3}
\ee
and
\be
F(G(y))=y\ .
\label{F4}
\ee
\end{theor}

\noindent{\bf Proof :} One easily proves by bounding the $w$ and $y$
factors
in $A_i(\cT,(i_l)_{l\in\cT})$ by their moduli and summing
the indices, starting with the leaves and progressing towards the root,
that
\be
|\cA_i(\cT)|\le {||w||_{\infty,1}}^V {||y||_{\infty}}^N
\ee
for any fixed tree $\cT$.
Therefore
\bea
\lefteqn{
\sum_{V,N\ge 0}
\frac{1}{V!N!}\sum_{\cT}
|\cA_i(\cT)|} & & \nonumber\\
 & & \le \sum_{V\ge 0}
\frac{(dV)!{||w||_{\infty,1}}^V {||y||_{\infty}}^{(d-1)V+1}}
{V!((d-1)V+1)!d!^V}\\
 & & \le ||y||_{\infty}
\sum_{V\ge 0}
\frac{(dV)!}{V!((d-1)V)!}
{\lp\frac{||w||_{\infty,1}{||y||_{\infty}}^{d-1}}{d!}
\rp}^V
\eea
and one simply uses
\be
\frac{(dV)!}{V!((d-1)V)!}\le
2^{dV}
\ee
to conclude the convergence proof and obtain the bound (\ref{F3}).
Now observe that, on the convergence domain
\be
G_i(y)=y_i+
\sum_{{V,N,\cT}\atop{V\ge 1}}
\frac{1}{V!N!}\cA_i(\cT)
\ee
where the last sum is over trees with at least one $w$-vertex linked
directly to the root.
This sum can be performed in the following way.
One chooses, among the $V$ internal $w$-vertices, the vertex $w_0\in
E_{\rm internal}$ which hooks to the root. This costs a factor $V$.
Then one divides the remaining vertices into an unordered collection of sets
$E_1,\ldots,E_d$ such that $E_i$ has $V_i$ $w$-vertices and
$N_i$ $y$-vertices.
This costs a factor
\[
\frac{1}{d!}
\frac{(V-1)!}{V_1!\ldots V_d!}
\frac{N!}{N_1!\ldots N_d!}\ .
\]
Finally one sums over all possible trees $\cT_1,\ldots\cT_d$
on $E_1\cup\{w_0\},\ldots,E_d\cup\{w_0\}$ as before.
The corresponding amplitudes do not depend on the location
of the sets $E_i$ in $E$, but only on the cardinalities $V_i$ and $N_i$.
Therefore
\bea
\lefteqn{
G_i(y)=y_i+\sum_{{V_1,\ldots,V_d\ge 0}\atop{N_1,\ldots,N_d\ge 1}}
\sum_{i_1,\ldots,i_d=1}^n
\sum_{\cT_1,\ldots,\cT_d}\frac{1}{V!N!}
} & & \nonumber\\
 & & \frac{V.(V-1)!.N!}{d!V_1!\ldots V_d!N_1!\ldots N_d!}
w_{i,i_1\ldots i_d}
\cA_{i_1}(\cT_1)\ldots\cA_{i_d}(\cT_d)\\
 &  & = y_i+\sum_{i_1,\ldots,i_d=1}^n
\frac{1}{d!}w_{i,i_1\ldots i_d}
G_{i_1}(y)\ldots G_{i_d}(y)\\
 &  & = y_i+H_i(G(y))
\eea
from which (\ref{F4}) follows.
\endproof

As a result the Taylor series of $G$ at the origin
is the right compositional inverse of $F$. Now algebraic combinatorialists
might not be too impressed by this since one can readily rewrite formula
(\ref{F5}) under the form given by Towber~\cite{Wright1} or
Singer~\cite{Singer}.
So the series expansion of the formal inverse itself is not new.
To obtain a real improvement on previous approaches one has to return
to the more fundamental equation (\ref{F1}) and really consider the
``integrals'' appearing in it as, well, {\em integrals}
on which one can try all the tools of ordinary calculus:
integration by parts, change of variables$\ldots$
For an example of the mathematical utility of this way of proceeding,
see~\cite{BarNatan1,BarNatan2,Rozansky}.

\noindent{\bf Remark :}
Note that the generalized forest formula of Towber~\cite{Wright1}
can be easily derived from the perturbation expansion of higher
correlation functions
$<\ph_1^{\al_1}\ldots\ph_n^{\al_n}>$,
where we used the standard statistical mechanics notation
\be
<\Om({\Br\ph},\ph)>
\eqdef
\frac{1}{Z}
\int
d{\Br\ph}d\ph\ \Om({\Br\ph},\ph)
e^{-{\Br\ph}\ph+\frac{1}{d!}{\Br\ph}w\ph^d
+{\Br\ph}y}\ \ .
\ee

\section{Comments on the Jacobian conjecture}

\subsection{What does the constant Jacobian condition mean?}

Suppose that
\be
F_i(x)=x_i-\sum_{j_1,\ldots,j_d=1}^n
\frac{1}{d!}
w_{i,j_1\ldots j_d}
x_{j_1}\ldots x_{j_d}
\ee
is such that $JF(x)=1$ for all $x$.
Several conclusions can be drawn from this constraint.
One that is due to V. Rivasseau is that our QFT model
is {\em self-normalized}.
In other words
\be
Z=1\ .
\ee
Indeed, by writing the Feynman diagram expansion of
\be
Z=\int
d{\Br\ph}d\ph\ 
e^{-{\Br\ph}\ph+\frac{1}{d!}{\Br\ph}w\ph^d
+{\Br\ph}y}
\ee
one can easily show that
\be
\log Z
=\sum_{k\ge 1}\frac{1}{k}
tr \left[
M(G(y))^k
\right]
\ee
where $M(x)$ is the matrix with entries
\bea
M_{ij}(x) & \eqdef &
\frac{\partial H_i}{\partial x_j}(x)\\
 & = & \frac{1}{(d-1)!}\sum_{j_1,\ldots,j_{d-1}=1}^n
w_{i,j j_1\ldots j_{d-1}}
x_{j_1}\ldots x_{j_{d-1}}
\eea
that is
\bea
Z & = & \exp \lp -tr\ \log \lp I-M(G(y))\rp\rp \\
 & = & \frac{1}{det\left[\frac{\partial F}{\partial x}(G(y))
\right]}\\
 & = & \frac{1}{JF(G(y))}\ .
\eea
In the case where $F_i(x)=x_i-H_i(x)$ with
the $H_i(x)$
homogenous of the same degree $d$, it is easy to show that
the Jacobian condition is equivalent to
$M(x)$ being nilpotent for all $x$ (see~\cite{BassCW}).
There are essentially two ways to express this
\be
M(x)^n=0,\ \ \ \forall x\in \CC^n
\label{F6}
\ee
or
\be
tr\lp M(x)^k\rp=0,\ \ \ \forall k\ge 1,
\forall x\in\CC^n
\label{F7}
\ee

Equation (\ref{F6}) means that when one considers a {\em chain}
(or caterpillar) diagram like
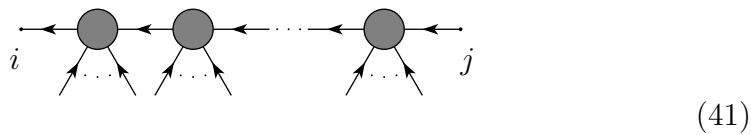
\be
\figput{chaindess}
\label{chain}
\ee
with $n$ $w$-vertices, its contribution, for fixed $i$ and $j$,
is zero after {\em symmetrization} of the indices of the
$n(d-1)$ incoming lower legs.

Equation (\ref{F7}) means that {\em loop} diagrams like
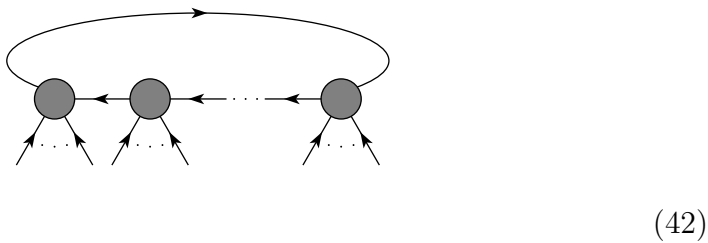
\be
\figput{loopdess}
\label{loop}
\ee
with $k\ge 1$ $w$-vertices, vanish after {\em symmetrization}
of the indices of the $k(d-1)$ incoming legs.
The formal inverse approach to the Jacobian conjecture can now be rephrased
as the following

\medskip
\noindent{\bf Problem :}
Show {\em explicitly} that in the polynomial algebra $\CC[w]$
with indeterminates given by the tensor elements $w_{i,j_1\ldots j_d}$,
the $y=0$ connected correlation functions
\[
<\ph_i{\Br\ph}_{j_1}\ldots{\Br\ph}_{j_N}>_{y=0}^c
\]
belong to the {\em radical} of the ideal generated
by the {\em symmetrized} chains and/or loops, provided the degree $N$
is large enough.

\medskip
This statement is by the Hilbert nullstellensatz
{\em equivalent} to the Jacobian conjecture.
It is even a theorem due to S. Wang~\cite{Wang} in the ($d=2$)
quadratic case.
The proof is non constructive however, and an explicit combinatorial
argument is an urgent desideratum.

\subsection{Chains and/or loops?}

Let $\Gc$ be the ideal of $\CC[w]$ generated by the symmetrized chains of
length $n$, and let $\Gl$ be the ideal generated by the symmetrized loops
of length $k\ge 1$.
While it is very tempting to work with $\Gc$, it seems more fundamental to
use $\Gl$.
This conclusion is implicit in~\cite{Wright1}. Indeed the author uses the
diagonal minor sums, i.e. the elementary symmetric functions of the
eigenvalues, to express the nilpotence of $M(x)$, instead of
the matrix elements of $M(x)^n$.
We use the loops, that is the Newton power sums of the eigenvalues,
which makes no difference since our ground ring is $\CC$.
Note that $\Gc\subset\Gl$: this is the Cayley-Hamilton theorem,
i.e. ``the Jacobian problem for $d=1$''!
But we also have $\Gl\subset\sqrt{\Gc}$,
trivially because a nilpotent matrix
must have zero eigenvalues and therefore the Newton sums of these eigenvalues
are zero.
It is very instructive to understand these two elementary statements
in a purely combinatorial way. Regarding the first inclusion,
we were surprised to find in the recent literature a combinatorial proof,
with a flavor of loop-erased random walk, of the eminently classical
Cayley-Hamilton theorem~\cite{Straubing}.
As for the second inclusion, there is a very nice explicit Fermionic
proof~\cite{CalanM} that, for a generic $n\times n$ matrix
$N$, $(tr\ N)^{k(n-1)+1}$ is in the ideal generated by the matrix
elements of $N^k$.

To see why the ideal $\Gc$ is tempting to work with, we need to recall a
theorem, first conjectured by Wang in the quadratic case~\cite{Wang},
and proved in full generality by O. Gabber (see~\cite{BassCW}).

\begin{theor}
If $F:\CC^n\rightarrow \CC^n$ is globally invertible with
polynomial inverse $G$ then
\be
{\rm deg}\ G\le({\rm deg}\ F)^{n-1}
\ee
where ${\rm deg}\ F\eqdef\max_{1\le i\le n}({\rm deg}\ F_i)$ and
likewise for $G$.
\end{theor}

In our context, this means that the vanishing of the connected
correlation functions
$<\ph_i{\Br\ph}_{j_1}\ldots{\Br\ph}_{j_N}>_{y=0}^c$
should happen as soon as $N>d^{n-1}$. Note that this bound is saturated by
the well known triangular example given by
$F_i(x)=x_i-x_{i+1}^d$, for $1\le i\le n-1$, and $F_n(x)=x_n$.
But $d^{n-1}$ is the maximal number of leaves of those of our trees
which have a depth less than or equal to $n-1$.
If the chains in (\ref{chain})
needed not be symmetrized, the Jacobian conjecture would be trivial!
Indeed, a tree with more than $d^{n-1}$ leaves must
have a chain
of length at least $n$, going from the root $\ph_i$ to one
of the leaves ${\Br\ph}_{j_\al}$.
This observation, which goes back to~\cite{BassCW},
was likely the main impetus
behind the formal inverse approach.

Remark that if we condition the sum over Feynman diagrams for
the correlations $<\ph_i{\Br\ph}_{j_1}\ldots{\Br\ph}_{j_N}>_{y=0}^c$,
by requiring the path between the root $\ph_i$
and a {\em specified} leaf ${\Br\ph}_{j_\al}$ to be of a certain length
$\ge n$; the branches will be automatically symmetrized and the
result would be zero. The problem is that we cannot know {\em in advance}
which
leaf will be linked to the root by a {\em long chain}.

In relation to previously used formal inversion formulas, let us mention
that it is against QFT wisdom to mix the index space $\{1,\ldots,n\}$
and the abstract space $E$ that labels the vertices, as far as the
combinatorics are concerned. From a QFT point of view, which admittedly
is only one
among many on the Jacobian conjecture, it is unnatural to use sums
over colored or planar objects, as this reduces symmetry in the
resulting expansion instead of enhancing it. We nevertheless concede
the point that planarity can serve to ``locate'' the long chain,
and order the trees accordingly, which is the main ingredient of the
combinatorial ``tour de force'' of~\cite{Singer}.

One of the cases treated in the latter article is that of
$F_i(x)=x_i-H_i(x)$, with the $H_i$ homogenous of the same degree $d$
and the matrix $M(x)$ nilpotent of order $2$.
This has already been treated in~\cite{BassCW}
and~\cite{ChengSW} for instance, but
let us sketch how to prove this result with our QFT model.
The argument is adapted from an idea by V. Rivasseau.

First perform the translation change of variables
$\ph\rightarrow \ph+y$, ${\Br\ph}\rightarrow{\Br\ph}$ in (\ref{F1})
to get, using $Z=1$,
\be
G_i(y)=y_i+
\int
d\mu({\Br\ph},\ph)\ \ph_i
e^{{\Br\ph}H(\ph+y)}
\label{F8}
\ee
where $d\mu({\Br\ph},\ph)\eqdef d{\Br\ph}d\ph\ e^{-{\Br\ph}\ph}$.
This unorthodox change of variables used in (\ref{F8}),
which treats $\ph$ and $\Br\ph$ as independent variables and not as
complex conjugates of one another, 
can be justified {\em a posteriori} by comparing the
diagrammatic expansions on both sides of the equation.
One can integrate the source $\ph_i$ by parts to get
\bea
G_i(y) & = & y_i+\int
d\mu({\Br\ph},\ph)\ \frac{\partial}{\partial{\Br\ph}_i}
e^{{\Br\ph}H(\ph+y)} \\
 & = & y_i+\int
d\mu({\Br\ph},\ph)\ H_i(\ph+y)
e^{{\Br\ph}H(\ph+y)} \\
 & = & y_i+ {\left. \int
d\mu({\Br\ph},\ph)\ H_i(\ph+y)
e^{s{\Br\ph}H(\ph+y)}\right|}_{s=1}\ \ .
\eea
Then, interpolate between $s=1$ and $s=0$ to get
\be
G_i(y)= y_i+\int
d\mu({\Br\ph},\ph)\ H_i(\ph+y)
+\int_0^1 ds\ \Om_i(s,y)
\label{F9}
\ee
where
\be
\Om_i(s,y)\eqdef \int
d\mu({\Br\ph},\ph)\ H_i(\ph+y)
\left[
{\Br\ph}H(\ph+y)
\right]
e^{s{\Br\ph}H(\ph+y)}\ \ .
\ee
Notice that the second term of (\ref{F9})
reduces to $H_i(y)$, whereas for the third we have, by integrating the
$\Br\ph$ by parts
\bea
\Om_i(s,y) & = &
\sum_{j=1}^n
\int
d\mu({\Br\ph},\ph)\ \frac{\partial}{\partial\ph_j}
\lp
H_i(\ph+y)
H_j(\ph+y)
e^{s{\Br\ph}H(\ph+y)}
\rp\\
 & = & \Om_i^1(s,y)+
\Om_i^2(s,y)+
\Om_i^3(s,y)
\eea
where
\be
\Om_i^1(s,y)\eqdef
\int
d\mu({\Br\ph},\ph) \lp
\sum_{j=1}^n
M_{ij}(\ph+y)H_j(\ph+y)
\rp
e^{s{\Br\ph}H(\ph+y)}
\ee
\be
\Om_i^2(s,y)\eqdef
\int
d\mu({\Br\ph},\ph)\ H_i(\ph+y)\lp
\sum_{j=1}^n
M_{jj}(\ph+y)
\rp
e^{s{\Br\ph}H(\ph+y)}
\ee
and
\be
\Om_i^3(s,y)\eqdef
\int
d\mu({\Br\ph},\ph)\ H_i(\ph+y)\lp
\sum_{j,k=1}^n
H_j(\ph+y)s{\Br\ph}_k
M_{kj}(\ph+y)
\rp
e^{s{\Br\ph}H(\ph+y)}
\ee
Note that $\Om_i^2(s,y)=0$ since it contains
$\sum_{j=1}^n M_{jj}(x)=tr\ M(x)$
at $x=\ph+y$, and
$M(x)$ is nilpotent.
Likewise, $\Om_i^1(s,y)$ and $\Om_i^3(s,y)$ vanish
since the integrand contains a factor of the form
\bea
\sum_{j=1}^n
M_{kj}(x)H_j(x) & = & \frac{1}{d}\sum_{j,l=1}^n
M_{kj}(x)M_{jl}(x)x_l\\
 & = & \frac{1}{d}\sum_{l=1}^n
{\left[
M(x)^2
\right]}_{kl}x_l\\
 & = & 0
\eea
where we used Euler's identity for the homogenous
$H_i$'s, and the fact that $M(x)$ is nilpotent of order 2.
As a result $G_i(y)=y_i+H_{i}(y)$.

\subsection{The Pauli exclusion principle}

In order to be able to prove the Jacobian conjecture by purely combinatorial
means,
one needs to exhibit a volume effect similar to the Pauli
exclusion principle, as otherwise one would not see
the finiteness of the index set $\{1,\ldots,n\}$ within
the strictly tensorial Feynman diagrammatic notation were indices
are contracted i.e. summed over.
One would love to have Fermions, instead Bosons, entering the picture.
Let us mention three, typically field theoretic, ideas 
that have not been pursued in previous attempts with the formal inverse
approach, and which deserve further investigation.

\subsubsection{Supersymmetry}

One way to introduce Fermions in a purely Bosonic
model is to exhibit a supersymmetry. 
If this could be done; it would probably be the ``voie royale''
towards understanding the conjecture.
Unfortunately we have not been able to make much headway 
in this direction so far.
Let us simply mention a strange feature of our model
that hints towards a hidden supersymmetry.
As a result of our choice of vertices and the fact that the
propagators are directed,  the perturbation expansion of the
one-point function is reduced to a tree graph expansion.
This means that the semi-classical expansion of our model
around the ``false vacuum'' $\ph=0$, ${\Br\ph}=0$ is {\em exact}.
Besides, the ``integrals'' in  (\ref{F1})  
which are supposed to be over $\CC^n$ reduce to the contribution of
a single critical point: the ``true vaccum'' obtained by solving
\be
\frac{\partial}{\partial \ph_i}
\lp
{\Br\ph}F(\ph)-{\Br\ph}y
\rp
=0
\ee
and
\be
\frac{\partial}{\partial{\Br\ph}_i}
\lp
{\Br\ph}F(\ph)-{\Br\ph}y
\rp
=0
\ee
that is $\ph=G(y)$ and ${\Br\ph}=0$.
This is reminiscent of the Duistermaat-Heckman theorem~\cite{DuistermaatH}
which is known to involve supersymmetry (see~\cite{Witten}).

\subsubsection{Renormalization}

The Gabber inverse degree bound, together with the previously given
example that saturates it, suggest that the sought exclusion principle
has to act along the chains from the root to the leaves of the trees
but not across, i.e. within generations.
This is quite odd in view of the eventual introduction of Fermionic
variables in our model.
This however hints to the possibility that the problem
may come from ``divergent'' two-point subgraphs i.e.
parts of the diagrams that look like
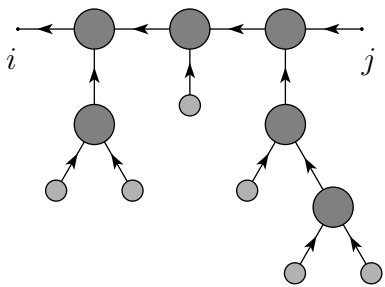
\be
\figput{renorm}
\ee
where the to indices $i$ and $j$ {\em coincide}.
This leads to the following.

\medskip
\noindent{\bf Question :} Is there a way to eliminate these ``divergent''
pieces by adding, to the ``action'' ${\Br\ph}F(\ph)-{\Br\ph}y$,
counter-terms that are made of symmetrized loops?

\medskip
This is possible for $d=1$, that is the Cayley-Hamilton theorem,
but does not seem to be the consequence of a natural
renormalization condition on the two-point correlation function.
It would be interesting to explore this idea using the new point of
view on renormalization pioneered by
A. Connes and D. Kreimer~\cite{ConnesK}
since one of their motivations was the study of formal
diffeomorphisms which is clearly related to our subject material.

\subsubsection{Reverse Mayer expansion}

We have repeatedly mentioned the Cayley-Hamilton theorem as the $d=1$
case of the Jacobian conjecture.
Let us now make this more precise.
If $d=1$ then $F_i(x)=x_i-H_i(x)$ with
$H_i(x)=\sum_{j=1}^n w_{i,j}x_j$ and $M(x)=M=(w_{i,j})_{1\le i,j\le n}$.
When $y=0$, the only interesting connected correlation functions
$<\ph_i{\Br\ph}_{j_1}\ldots{\Br\ph}_{j_N}>_{y=0}^c$ are for $N=1$, and
their Feynman diagram expansion only generates chains of the form
\be
\figput{matchain}
\ee
If we use the loop ideal $\Gl$ in the formulation of the problem
in section III.1, then the ensuing statement that chains of length $\ge n$
belong to $\sqrt{\Gl}$, or simply $\Gl$ here, {\em is}
the Cayley-Hamilton theorem. Therefore any combinatorial solution of this
problem, for general $d$, should {\em at the very least}
reproduce this well known result.
Now the Cayley-Hamilton theorem is the statement that the derivative
of
$\frac{1}{Z}=det(I-M)$, with respect to a matrix element of $M$,
has a vanishing component in degree $n$ in the $M_{ij}$ variables.
This in turn stems from the fact that $\frac{1}{Z}$ is itself
a {\em polynomial}
of
degree $n$ in $M$.
We all know this from our ``Fermionic/determinantal'' upbringing,
but let us suppose for a moment that all we know is Bosons
and the only expression available to us is
\be
Z=\int d{\Br\ph}d\ph\ e^{-{\Br\ph}\ph+{\Br\ph}M\ph}\ \ .
\ee
It is amusing to prove that $\frac{1}{Z}$ is a polynomial in $M$
of
degree at most $n$, using this formula.
Let us explain an answer that is inspired from the Mayer
expansion in statistical
mechanics (see~\cite{Ruelle}), and constructive QFT
(see~\cite{Rivasseau}).

Consider the partition function $\cZ$ of a gas of propagators
\be
\figplace{propag}{0 in}{-0.07 in}=\de_{ij}
\ee
and vertices
\be
\figplace{matvertex}{0 in}{-0.45 in}=M_{ij}
\ee
with two kinds of interactions.

- The propagators attach to the vertices in all possible ways, as
when applying Wick's theorem.

- The propagators can interact by {\em Mayer-links}
represented by a squiggly line
\be
\figput{mayerlink}
\ee
and carrying a factor $-\de_{ij}$
where $i$ and $j$ are the indices of the two propagators.
One also assigns by hand a factor $-1$ per propagator-vertex loop,
which we confess is cheating a bit.
$\cZ$ is therefore a sum of objects like
\be
\figput{mayergraph}
\label{Fig3}
\ee
with the appropriate symmetry factors, that is $\frac{1}{k}$
for each oriented loop with $k$ vertices,
an overall $\frac{1}{m!}$ if there are $m$ loops, and a global sign
$(-1)^m$.
On the one hand, if one sums over the structure of Mayer-links, with
fixed contraction
of the propagators to the vertices,
one rebuilds a ``hardcore constraint'' factor $(1-\de_{ij})$
{\em for each pair of propagators}, which is the opposite of the operation
one would do in a standard Mayer expansion.
Since the available index space $\{1,\ldots,n\}$
has cardinality $n$, there cannot be more than
$n$ propagators in a nonvanishing graph,
ergo $\cZ$ is of degree at most $n$ in
$M$.
On the other hand, $\log \cZ$, is a sum over connected objects
like those of (\ref{Fig3}), where connectedness involves
{\em both} types of lines. When one sums over Wick contractions, with
fixed configuration of Mayer-links, the result
is zero as soon as there is at least one Mayer-link, because
of the following {\em exchange move}
along the Mayer-link.
\be
\figplace{exchangeleft}{0 in}{-0.45 in}
\;\;\Longrightarrow\;\;
\figplace{exchangeright}{0 in}{-0.45 in}
\ee
Indeed such a move does not affect the amplitudes, but
modifies the loop count by one unit and thus the sign of the graph
(compare with~\cite{Straubing}).
As a result, $\log \cZ$ is a sum over single loops without Mayer-links and
with a single (-1) factor.
That is $\log \cZ = -\log Z$ which concludes the argument.

An interesting question raised by this approach is

\medskip
\noindent{\bf Question : }
Is there a hyperdeterminant, in the sense of~\cite{GelfandKZ}, that would
play, when $d\ge 2$, the role played by $\frac{1}{Z}$
when $d=1$, and that would, upon derivation with respect to a tensor
element $w_{i,j_1\ldots j_d}$ around a solution of the Jacobian condition,
give some finiteness information on our correlation functions?

\subsubsection{Are these three ideas different, really?}

Although we have no precise unified framework to propose at the moment,
we are tempted to say no.
The interplay between supersymmetry, renormalization and
the Mayer expansion is quite mysterious and is probably related to
the combinatorics of the symmetric group and the inclusion-exclusion
principle.

We will conclude by pointing out a few references where some clues on
these relationships might be found.
Mayer expansions involve coefficients which are M{\oe}bius
functions of certain partition lattices (see~\cite{Rota}).
These coefficients can be calculated by an analog of the classical
forest formula of
Zimmermann in renormalization theory
(see the introduction to chap. 4 of~\cite{Abdesselam3}).
They can also be expressed using the so called
Brydges-Kennedy forest formula~\cite{BrydgesK} that was
first proved there using
the Hamilton-Jacobi equation.
In~\cite{AbdesselamR1} we gave a purely algebraic proof of the latter
using some partial fraction combinatorial identities
(Lemma II.2 in~\cite{AbdesselamR1}).
Such identities have been given a very elegant interpretation
in terms of minimal factorizations of permutations as a product
of transpositions~\cite{Lafforgue}.
The global sign in the Mayer coefficients
is $(-1)^k$ where $k$ is the number of edges in the forest.
This sign obviously becomes
the signature of the permutation in the latter interpretation.
This strongly
suggests a relationship between Mayer expansions and
Fermions which was also alluded to in~\cite{Rivasseau}.
Note finally that the Brydges-Kennedy identity was considerably generalized
in~\cite{AbdesselamR2} (section III.2.1), where
critical use is made of shuffles and a kind of Chen's lemma
(see remarks following the proof of Lemma 9), although we
did not know this at the time.

\section{Conclusion}

We hope to have provided enough evidence that the Jacobian conjecture
is a very beautiful combinatorial challange, where mathematicians, either
conceptually of computationally inclined, and theoretical physicists
could fruitfully share their knowledge.
While future progress on the conjecture itself is still uncertain,
there are bound to be benefits from such an interdisciplinary
collaboration on this problem.

\end{document}